\documentclass[twocolumn]{autart}    

\usepackage{graphicx}          

\usepackage{graphicx}
\usepackage{amssymb}
\usepackage{amsmath}
\usepackage[T1]{fontenc}
\usepackage{subfigure}
\usepackage{epsfig}
\usepackage{color}
\usepackage{enumitem}
\usepackage{float}
\usepackage{soul}
\usepackage{wrapfig}
\usepackage{arydshln}
\usepackage[justification=centering]{caption}
 
 \usepackage{wrapfig}
    \usepackage{graphicx}

\usepackage{url}

\newtheorem{Remark 1}{Remark}
\newtheorem{Remark 2}[Remark 1]{Remark}
\newtheorem{Remark 3}[Remark 1]{Remark}
\newtheorem{Remark 4}[Remark 1]{Remark}
\newtheorem{Remark 5}[Remark 1]{Remark}
\newtheorem{Remark 6}[Remark 1]{Remark}
\newtheorem{Remark 7}[Remark 1]{Remark}

\newtheorem{Lemma 1}{Lemma}
\newtheorem{Lemma 2}[Lemma 1]{Lemma}
\newtheorem{Lemma 3}[Lemma 1]{Lemma}
\newtheorem{Lemma 4}[Lemma 1]{Lemma}
\newtheorem{Lemma 5}[Lemma 1]{Lemma}
\newtheorem{Lemma 6}[Lemma 1]{Lemma}
\newtheorem{Lemma 7}[Lemma 1]{Lemma}
\newtheorem{Assumption 1}{Assumption}
\newtheorem{Assumption 2}[Assumption 1]{Assumption}
\newtheorem{Assumption 3}[Assumption 1]{Assumption}
\newtheorem{Assumption 4}[Assumption 1]{Assumption}
\newtheorem{Definition 1}{Definition}
\newtheorem{Theorem 1}{Theorem}
\newtheorem{Theorem 2}[Theorem 1]{Theorem}
\newtheorem{Theorem 3}[Theorem 1]{Theorem}
\newtheorem{Theorem 4}[Theorem 1]{Theorem}
\newtheorem{Theorem 5}[Theorem 1]{Theorem}
\newtheorem{Theorem 6}[Theorem 1]{Theorem}
\newtheorem{Theorem 7}[Theorem 1]{Theorem}
\newtheorem{Theorem 8}[Theorem 1]{Theorem}
\newtheorem{Theorem 9}[Theorem 1]{Theorem}
\newtheorem{Theorem 10}[Theorem 1]{Theorem}

\begin{document}

\begin{frontmatter}

\title{Decentralized Nonconvex Optimization with Guaranteed Privacy and
Accuracy\thanksref{footnoteinfo}} 

\thanks[footnoteinfo]{The work was supported in part by the National Science Foundation under
Grants ECCS-1912702, CCF-2106293, CCF-2215088, and CNS-2219487. \\This paper was not presented at any IFAC
meeting. Corresponding author Y.~Q.~Wang. Tel. +1-864.656.5923.}

\author[YQ]{Yongqiang Wang}\ead{yongqiw@clemson.edu},    
\author[TB]{Tamer Ba\c{s}ar}\ead{basar1@illinois.edu}

\address[YQ]{Department of Electrical and Computer Engineering, Clemson University, Clemson, SC 29634, USA}  
\address[TB]{Coordinated Science Lab, University of Illinois
at Urbana-Champaign, Urbana, IL 61801, USA}             


\begin{abstract}                          
Privacy protection  and nonconvexity are two challenging problems in
decentralized optimization and learning involving sensitive data.
Despite some recent advances  addressing each of the two problems
separately, no results have been reported that  have theoretical
guarantees on both privacy protection and saddle/maximum avoidance
in decentralized nonconvex optimization. We propose a new algorithm
for decentralized nonconvex optimization that can enable both
rigorous differential privacy and   saddle/maximum avoiding
performance. The new algorithm allows the incorporation of
persistent additive noise to enable rigorous differential privacy
for  data samples, gradients, and intermediate optimization
variables without losing provable convergence, and thus
circumventing the dilemma of trading accuracy for privacy in
differential privacy design. More interestingly, the
 algorithm is  theoretically proven to be able to  efficiently { guarantee accuracy by avoiding}  convergence to local
maxima and saddle points, which has not been reported before in the
literature on decentralized nonconvex optimization. The algorithm is
efficient in both communication (it only shares one variable in each
iteration) and computation (it is encryption-free), and hence is
promising for large-scale nonconvex optimization and learning
 involving high-dimensional optimization parameters. Numerical experiments
 for
both a decentralized estimation problem and an Independent Component
Analysis (ICA) problem
 confirm the effectiveness of the proposed approach.
\end{abstract}

\end{frontmatter}

\section{Introduction}
Decentralized optimization is gaining increased traction across
disciplines due to its fundamental role in cooperative control
\cite{yang2019survey}, distributed sensing
\cite{bazerque2009distributed}, multi-agent systems
\cite{raffard2004distributed}, sensor networks
\cite{zhang2017distributed}, and large-scale machine learning
\cite{tsianos2012consensus}. In many of these applications, the
problem can be formulated in the following general form, in which a
network of $m$ agents cooperatively solve a common optimization
problem through on-node computation and local communication:
\begin{equation}\label{eq:optimization_formulation1}
\min\limits_{\theta\in\mathbb{R}^d} F(\theta)\triangleq
\frac{1}{m}\sum_{i=1}^m f_i(\theta)
\end{equation}
where the optimization variable $\theta$ is  common to all agents
but each $f_i(\theta):\mathbb{R}^d\rightarrow\mathbb{R}$ is a local
objective function private to agent $i$.

Plenty of results have been reported to solve the above
decentralized optimization problem since the seminal work of
\cite{tsitsiklis1984problems}, with some of the popular approaches
including gradient-descent (e.g.,
\cite{nedic2009distributed,qu2017harnessing,wang2022gradient}),
distributed alternating direction method of multipliers (e.g.,
\cite{shi2014linear,zhang2019admm}),  and distributed Newton methods
(e.g., \cite{wei2013distributed}). Results have also emerged
incorporating various communication and computation constraints in
decentralized optimization (e.g.,
\cite{lin2016distributed,cao2020decentralized,cao2021decentralized}).
 Most of the reported results focus on convex
objective functions, whereas results are relatively sparse for
nonconvex objective functions. However,  in many practical
applications, the objective functions are essentially nonconvex. For
example, in the resource allocation problem  of communication
networks, the utility functions cannot be modeled by convex/concave
functions when the communication traffic is
  non-elastic  \cite{tychogiorgos2013non}; in most machine  learning applications, the objective functions are
    essentially  nonconvex due to the presence of multi-layer
    neural networks \cite{tsianos2012consensus}; in policy optimization for linear-quadratic regulator
    \cite{fazel2018global}
    as well as for robust and risk-sensitive control \cite{zhang2020policy}, nonconvex optimization naturally arises.
  Therefore, it is imperative to study decentralized optimization under nonconvex
  objective functions.

  In recent years, results have emerged on decentralized nonconvex optimization
  \cite{bianchi2012convergence,di2016next,tatarenko2017non,wai2017decentralized,zeng2018nonconvex}, which address
   the  convergence of participating agents' optimization variables to
  a first-order stationary point of the global objective function.
  Nevertheless, these results do not address the avoidance of saddle
  points (stationary points that are  not local extrema), which is
 a major concern in many nonconvex optimization problems \cite{ge2015escaping}. For example,
    in machine learning applications, the main bottleneck in parameter optimization is not due to the existence of multiple local minima, but the
    existence of many saddle points which  trap gradient updates
    \cite{ge2015escaping}. To escape from saddle points, classical
    approaches resort to second-order  information, in particular the Hessian matrix of second derivatives (see, e.g., \cite{nesterov2006cubic,curtis2017trust}).
    The Hessian matrix based approach,  however, incurs a high cost in both computation and storage in every iteration since the Hessian matrix scales quadratically with
    the
    dimension of optimization variables, which can be hundreds of millions  in modern deep learning applications   \cite{tang2020communication}.
     Recently, first-order gradient methods have been shown to be able
     to escape saddle points with the help of random perturbations
     in centralized
     optimization (see, e.g., \cite{ge2015escaping,jin2021nonconvex}). However, it is unclear if this is still true in
     decentralized nonconvex optimization since the decentralized  architecture
     brings in
     fundamental differences in optimization dynamics. For example, the saddle points of
     individual objective functions $f_i(\cdot)$ in decentralized optimization are  different from those of
 the aggregated objective function $F(\cdot)$, which is the only  function that needs to be considered in centralized
 optimization. Furthermore, the inter-agent coupling also
 complicates the optimization dynamics.
      Note that random initialization has been shown to be able to asymptotically avoid saddles in centralized nonconvex optimization \cite{lee2016gradient},
      which is further extended to the decentralized case in \cite{daneshmand2020second}. However, the result in \cite{du2017gradient} shows that
      this approach to  avoiding saddles may take exponentially long time, rendering
      it
      impractical.

     In this paper, we propose an approach to avoiding maxima/saddles in decentralized nonconvex optimization
     by leveraging differential privacy design. More specifically, we propose a new algorithm for first-order decentralized
     optimization that enables exploiting differential-privacy noise to achieve guaranteed saddle-avoiding performance without losing provable convergence.
      This is significant because  differential-privacy noise
is generally known to sacrifice provable convergence while enabling
privacy
protection \cite{bagdasaryan2019differential}. 
Moreover, we rigorously establish that the differential-privacy
noise can prevent the algorithm from converging to undesired
stationary points such as saddle points within polylogarithmic time,
and hence can enhance optimization accuracy. This extends recent
results on saddle avoidance in centralized nonconvex optimization
\cite{jin2021nonconvex}, and to our knowledge,  has not been
reported for decentralized  optimization with guaranteed convergence.
It is
     worth noting that although {diminishing noises}, known as annealing,
     has been shown to be able to
     facilitate global convergence in {distributed nonconvex optimization \cite{swenson2019annealing,swenson2020distributed}},
such an approach may not be efficient as  it could result in
convergence time increasing exponentially with the dimension
     of optimization variables \cite{gayrard2004metastability}.

     The considered privacy protection aspect is becoming an increasingly pressing
     need in decentralized optimization involving sensitive
     data.
  For example, in sensor network based localization, the positions
  of  sensor agents should  be kept private in sensitive or hostile
  environments \cite{zhang2019admm,huang2015differentially}. This
  requires that
  decentralized-optimization based localization approaches
  protect the privacy of individual agents' gradients, which are
  linear functions of sensor positions and whose disclosure will directly reveal a
sensor's position \cite{zhang2019admm}.   Another example
underscoring the importance of privacy preservation in decentralized
optimization is machine learning, where the data involved  may
contain sensitive information such as medical records or salary
information \cite{yan2012distributed}. In fact, recent  results in
\cite{zhu2019deep} (as well as our own result \cite{wang21,wang2022decentralized}) show
that without a strong privacy mechanism in place, an adversary can
precisely recover the raw data used for training through shared
gradients (pixel-wise accurate for images and token-wise matching
for texts).

   The main contributions of this paper are as follows: 1)~We
   propose a new algorithm  for decentralized optimization that enables
   the
   achievement of  differential privacy  without losing provable convergence. This is
   significant
   since in general differential privacy has to trade algorithmic
   accuracy for privacy; 2) We rigorously establish that the proposed algorithm can { guarantee accuracy by efficiently avoiding} saddle
   points in decentralized nonconvex optimization under a diminishing
   stepsize, which, to our knowledge, has not been reported before.
   We would like to emphasize that allowing the stepsize to diminish
   with time
   is crucial to ensure   provable convergence  under  persistent differential-privacy noise; 3) We   prove
   that the proposed approach can  enable rigorous differential
   privacy for individual agents' data samples and shared gradients  under guaranteed convergence. Moreover,  the proposed approach
     also provides differential privacy for   participating agents' optimization
   variables. However, different from the differential-privacy protection on data samples and gradients, the differential-privacy protection on
optimization variables decreases with time and reduces to zero when
the algorithm converges. Note that this is completely acceptable
because the objective of decentralized optimization is to let
individual agents   learn the same optimal optimization variable.

The organization of the paper is as follows. Sec. II provides
formulation of  the problem. Sec. III presents the proposed
algorithm and Sec. IV provides a rigorous analysis on convergence,
including the guaranteed performance in avoiding maxima and saddle
points  in decentralized nonconvex optimization. Sec. V discusses
the privacy preservation performance of the algorithm. Sec. VI
presents numerical experimental results in both a decentralized
estimation application and an Independent Component Analysis (ICA)
application. Finally, Sec. VII concludes the paper.

{\bf Notations:} $\mathbb{R}^m$ denotes the Euclidean space of
dimension $m$. $I_d$ denotes the identity matrix of dimension $d$,
and ${\bf 1}_d$ denotes a $d$ dimensional column vector with all
entries equal to 1; in both cases  we suppress the dimension when
clear from the context. A vector is viewed as a column vector. For a
vector $x$, $x_i$ denotes its $i$th element.
$A^T$ denotes the transpose of matrix $A$. $\|x\|$ denotes the
standard Euclidean norm $\|x\|=\sqrt{\sum_{i=1}^d(x_i)^2}$ and
$\|x\|_1$ denotes the Taxicab norm $\|x\|_1=\sum_{i=1}^d|x_i|$. A
matrix is column-stochastic when its entries are nonnegative and
elements in every column add up to one. A square matrix $A$ is said
to be doubly-stochastic when both $A$ and $A^T$ are
column-stochastic.

\section{Problem Formulation}

\subsection{Objective functions}
In decentralized optimization, agent $i$'s objective function
$f_i(\cdot)$ is determined by its loss function and locally
accessible data samples. Therefore, we consider $f_i(\cdot)$ of the
following form
\begin{equation}\label{eq:f_i}
f_i(\theta)=\frac{1}{n_i}\sum_{j=1}^{n_i}\ell_i(\theta,s_{i,j})
\end{equation}
where $\ell_i(\cdot,\cdot)$ denotes the cost function of agent $i$,
$n_i$ denotes the number of data samples available to agent $i$, and
$s_{i,j}$ represents the $j$th data sample of agent $i$. We
represent the set of all data samples available to agent $i$ as
$\mathbb{D}_i$.

 We make the following
assumption on cost functions:
\begin{Assumption 2}\label{assumption:gradient}
Every $\ell_i(\cdot,\cdot)$ satisfies
$\lim\limits_{\|u\|\rightarrow\infty}\ell_i(u,\cdot)\rightarrow\infty$
and
has Lipschitz gradient and Lipschitz
  Hessian over $\mathbb{R}^d$, i.e., for some $\nu>0$ and $\rho>0$,
\[
\begin{aligned}
\|\nabla \ell_i(u,\cdot)-\nabla \ell_i(v,\cdot)\|\le \nu
\|u-v\|,\quad \forall u,v\in\mathbb{R}^d.
\\
\|\nabla \ell_i(\cdot,s_{p})-\nabla \ell_i(\cdot,s_{q})\|_1\le \nu
\|s_{p}-s_{q}\|_1,\quad \forall s_{p},s_{q}\in\mathbb{D}_i.
\\
\|\nabla^2 \ell_i(u,\cdot)-\nabla^2 \ell_i(v,\cdot)\|\le \rho
\|u-v\|,\quad \forall u,v\in\mathbb{R}^d.
\end{aligned}
\]
always hold for all $i$, where $\nabla^2 \ell_i(\cdot)$ denotes the
Hessian matrix of $\ell_i(\cdot,\cdot)$ with respect to the first
argument.
\end{Assumption 2}

From the definition of $f_i(\cdot)$ in (\ref{eq:f_i}), it can be
easily verified that $f_i(\cdot)$ always has Lipschtiz gradient and
Hessian, i.e.,
\[
\begin{aligned}
\|\nabla f_i(u,\cdot)-\nabla f_i(v,\cdot)\|\le \nu \|u-v\|,\quad
\forall u,v\in\mathbb{R}^d.\\
\|\nabla^2 f_i(u,\cdot)-\nabla^2 f_i(v,\cdot)\|\le \rho
\|u-v\|,\quad \forall u,v\in\mathbb{R}^d.
\end{aligned}
\]

{ The coercivity   assumption $\lim\limits_{\|u\|\rightarrow\infty}\ell_i(u,\cdot)\rightarrow\infty$ is used here because we need  the stochastic approximation theory in \cite{nevelson1976stochastic} to prove the avoidance of local maxima. It is also recently used in \cite{tatarenko2017non} to analyze the push-sum based distributed optimization under the assumption of no saddle points.   The Lipschitz   gradient and  Hessian condition is a standard assumption in saddle-avoidance studies \cite{ge2015escaping,jin2021nonconvex,daneshmand2018escaping}.

We  also assume that the gradient is bounded, which is commonly used in (distributed) nonconvex optimization  \cite{daneshmand2018escaping,moulines2011non,jiang2017collaborative,koloskova2019decentralized} and differential-privacy analysis \cite{huang2015differentially,abadi2016deep}:

\begin{Assumption 1}\label{As:bounded_gradient}
For every $i$, we always have
$
\|\nabla \ell_i(\cdot,\cdot)\|\leq G
$
for some positive constant $G<\infty$, which further implies $
\|\nabla f_i(\cdot)\|\leq G
$ according to the definition of $f_i(\cdot)$ in (\ref{eq:f_i}).
\end{Assumption 1}
\begin{Remark 1}
  Note that the Lipschitz function assumption in \cite{jiang2017collaborative} implies bounded gradients.
\end{Remark 1}

}
For a twice differentiable aggregated objective function $F(\cdot)$,
we call $\theta$ a stationary point if $\nabla F(\theta)=0$ holds. A
stationary point $\theta$ can be  viewed as belonging to one of
three categories:
\begin{itemize}
\item local minimum: there exists a $\gamma>0$ such that $F(\theta)\leq
F(\vartheta)$ for any $\|\vartheta-\theta\|\leq \gamma$;
\item local maximum: there exists a $\gamma>0$ such that $F(\theta)\geq
F(\vartheta)$ for any $\|\vartheta-\theta\|\leq \gamma$;
\item saddle point: neither of the above two cases is true, i.e., for any $\gamma>0$, there exist
$\vartheta_1$ and $\vartheta_2$ satisfying
$\|\vartheta_1-\theta\|\leq \gamma$ and $\|\vartheta_2-\theta\|\leq
\gamma$ such that $F(\vartheta_1)<F(\theta)<f(\vartheta_2)$.
\end{itemize}

Since distinguishing saddle points from local minima for smooth
functions is   NP-hard in general \cite{nesterov2000squared}, we
focus on a subclass of saddle points, i.e., strict saddle points:

\begin{Assumption 1}\label{assumption:strict_saddle}
All saddle points $\theta$ of the aggregated function $F(\cdot)$ are
strict saddles, i.e., the minimum (resp. maximum) eigenvalue of the
Hessian matrix $\nabla^2 F(\theta)$ at any saddle $\theta$  is
negative (resp. positive).
\end{Assumption 1}

A generic saddle point must satisfy that the minimum (resp. maximum)
eigenvalue of its Hessian matrix is non-positive (resp.
non-negative). Our assumption of strict saddles rules out the case
where the minimum or  maximum  eigenvalue of the Hessian matrix is
zero. A line of recent work in the machine learning literature shows
that for many popular models in machine learning, all saddle points
are indeed strict saddle points, with examples ranging from tensor
decomposition \cite{ge2015escaping}, dictionary learning
\cite{sun2016complete}, smooth semidefinite programs
\cite{boumal2016non}, to robust principal component analysis
\cite{ge2017no}.

Recently, \cite{ge2015escaping} and \cite{jin2021nonconvex} have
shown that in centralized nonconvex optimization, saddle points
could be avoided efficiently (in a polylogarithmic number of
iterations) by adding perturbations in the classical single-variable
gradient descent algorithm. In this paper, we extend this result to
the decentralized case and prove that the added differential-privacy
noise can be leveraged to avoid saddles without sacrificing
provable convergence. It is worth noting that the extension from
centralized optimization to the decentralized case is highly
nontrivial because the saddle points of the aggregated function
$F(\theta)$ are different from those of individual objective
functions $f_i(\theta)$. Furthermore, in decentralized optimization,
the interaction between agents brings in an additional element that
affects the evolution of dynamics around saddle points, which makes
state evolution analysis around saddle points more involved compared
with the centralized optimization case.

\subsection{Interaction topology}
We consider a network of $m$ agents. The agents interact on an
undirected graph, which can be described by a weight matrix
$W=\{w_{ij}\}$. More specifically, if agents  $i$ and  $j$ can
  interact with each other, then $w_{ij}$ is
positive. Otherwise, $w_{ij}$ will be zero. We assume that an agent
is always able to affect itself, i.e., $w_{ii}>0$ for all $1\leq
i\leq m$. The neighbor set $\mathbb{N}_i$ of agent $i$ is defined as
the set of agents $\{j|w_{ij}>0\}$. So the neighbor set of agent $i$
always includes itself. To ensure that the agents can cooperatively
solve the decentralized optimization problem
(\ref{eq:optimization_formulation1}), we make the following standard
assumption on the interaction topology:

\begin{Assumption 1}\label{assumption:W}
   $W=\{w_{ij}\}\in \mathbb{R}^{m\times m}$
   is symmetric and  satisfies ${\bf 1}^TW={\bf
  1}^T$, $W{\bf 1}={\bf
  1}$, and $\eta=\|W-\frac{{\bf 1}{\bf 1}^T}{m}\|<1$.
\end{Assumption 1}

The optimization problem (\ref{eq:optimization_formulation1}) can
now be reformulated as the following equivalent multi-agent
optimization problem:
\begin{equation}\label{eq:optimization_formulation2}
\min\limits_{x\in\mathbb{R}^{md}}f(x)\triangleq
\frac{1}{m}\sum_{i=1}^m f_i(x_i)\:\: {\rm s.t.}\:\:
x_1=x_2=\cdots=x_m
\end{equation}
where $x_i\in\mathbb{R}^d$ is the local estimate of agent $i$ about
the optimization solution and
$x=[x_1^T,x_2^T,\cdots,x_m^T]^T\in\mathbb{R}^{md}$ is the collection
of the estimates made by the agents.

\subsection{Privacy preservation in decentralized optimization}\label{se:privacy_formulation}

In decentralized optimization, the sensitive information that has to
be protected from disclosure could be the data samples (raw data),
gradients, or optimization variables. Although  data samples are not
directly shared in decentralized optimization, their information are
abstracted and embedded in gradients. For example, in
decentralized-optimization based rendezvous and localization,
disclosing the gradient  of an agent amounts to disclosing its
(initial) position \cite{zhang2019admm,huang2015differentially,wang2022tailoring},
which directly correlates with sampled range/angle measurements. In
machine learning, it has been shown that shared gradients can be
used by an adversary to reversely recover the raw data used for
training (pixel-wise accurate for images and token-wise matching for
texts) \cite{zhu2019deep,wang21,wang2022decentralized}. The optimization variables (models in
machine learning) could also carry sensitive information abstracted
from raw data. However, note that the objective of decentralized
optimization is for individual agents to learn the same optimal
optimization variable (model), and hence the   final consensual
optimization variable should be disclosed to all agents, and not be
a target of privacy protection.  Therefore, in this paper, we
restrict the privacy to individual agents' data samples and
gradients, and individual agents' intermediate optimization
variables (by ``intermediate,'' we mean the evolution of optimization
variables before achieving consensus among agents).

  We consider two potential attacks in decentralized optimization, which are the
two most commonly used models of attacks in privacy research
\cite{Goldreich_2}:
\begin{itemize}
\item \emph{Honest-but-curious attacks}  are attacks in which a participating
agent or multiple participating agents (colluding or not)
 follow all protocol steps correctly but are curious and
collect  all received intermediate data  to learn the sensitive
information about other participating agents.

\item \emph{Eavesdropping attacks}  are attacks in which an external eavesdropper eavesdrops upon all communication channels to
intercept exchanged messages so as to learn sensitive information
about the sending agents.
\end{itemize}

An honest-but-curious adversary (e.g., agent $i$) has access to the
internal state  $x_i$, which is unavailable to external
eavesdroppers. However, an eavesdropper has access to all shared
information in the network, whereas an honest-but-curious agent only
has   access to shared information that is destined to it.

In this paper, the proposed  new decentralized optimization
algorithm enables us to leverage differential-privacy noises to
facilitate the avoidance of maxima and saddle points in nonconvex
optimization. We adopt the popular definition of
($\epsilon,\delta$)-differential privacy following standard
conventions \cite{dwork2014algorithmic}:
\begin{Definition 1}
A randomized function $h(x)$  is $(\epsilon, \delta)$-differentially
private if for all $S\subset {\rm Range}(h)$ and for all $x,y$ with
$\|x-y\|_1\leq 1$, we   have
\[
{\rm Prob}(h(x)\in S)\leq e^\epsilon {\rm Prob}(h(y)\in S)+\delta
\]
where { ${\rm Range}(h)$ denotes the image (the set of all output values) of the function $h$ and} $\rm Prob(\cdot)$ denotes probability.
\end{Definition 1}

Note that $\epsilon$ and $\delta$ are always non-negative, and a
smaller $\epsilon$ (or $\delta$) corresponds to a stronger privacy
protection.

\section{The proposed decentralized optimization algorithm}

Conventional single-variable decentralized optimization algorithms
usually take the following form \cite{nedic2009distributed}:
\begin{equation}\label{eq:conventional_gradient_descent}
 x_i^{k+1}=\sum\nolimits_{j\in\mathbb{N}_i}w_{ij}x_j^k+\lambda^k g_i^k
\end{equation}
where $x_i^k$ denotes the local copy of optimization variable of
agent $i$ at iteration $k$, $\lambda^k$ is a positive scalar
denoting the stepsize, and
  $g^k_i$ denotes
the gradient of agent $i$ evaluated at $x_i^k$, i.e., $g_i^k=\nabla
f_{i}(x_i^k)$. It is well-known that under Assumption
\ref{assumption:gradient} and Assumption \ref{assumption:W}, when
$\lambda^k$ is such that  $\sum_{k=0}^{\infty}\lambda^k=\infty$ and
$\sum_{k=0}^{\infty}(\lambda^k)^2<\infty$, then all $x_i^k$ will
converge to the same optimal solution when $f(\cdot)$ is  convex.

However, in the above decentralized optimization  algorithm, agent
$i$ has to share $x^k_i$
 with all its neighbors $j\in\mathbb{N}_i$, which   breaches the privacy of optimization variable $x^k_i$.  Furthermore, if
an adversary has access to the optimization variable $x_i^k$ of
agent $i$
 and the updates that agent $i$ receives from all its neighbors
$x_j^k$ for $j\in\mathbb{N}_i$, then the adversary can easily infer
$g_i^k$   based on the update rule
(\ref{eq:conventional_gradient_descent}) and publicly known  $W$ and
$\lambda^k$. To protect the privacy of individual agents'
optimization variable $x^k_i$, existing decentralized optimization
approaches usually choose to inject additive noise on shared $x^k_i$
(see, e.g., \cite{huang2015differentially}), which, however, will
compromise the accuracy of the final optimization result.

We propose the following  decentralized optimization algorithm:
\begin{equation}\label{eq:algorithm_individual_agent}
 x_i^{k+1}=\sum_{j\in\mathbb{N}_i}w_{ij}(x_j^k-\lambda^kg_j^k)
\end{equation}

 The detailed implementation procedure for individual agents is
provided in Algorithm 1. {Compared with the conventional
decentralized optimization algorithm, it can be seen that instead of
letting agent $j$
 share $x_j^k$ with neighboring agents, we let agent $j$ share
$v_{ij}^k\triangleq w_{ij}(x_j^k-\lambda^k g_j^k)$ with all its
neighbors $i\in\mathbb{N}_j$. This new algorithm  has two advantages
over the conventional one in
(\ref{eq:conventional_gradient_descent}): First, it includes the
server based distributed optimization like federated learning as a
special case. More specifically, when all $x_j^k$ ($1\leq j\leq m$) are forced to be
the same, then different agents use the same parameter $x_j^k$ but
different local data sets to calculate gradients, the average of
which is used to update the universal state. This is exactly the
architecture used in federated learning. Note that since the conventional decentralized algorithms share $x_j^k$ among participating agents, they cannot be used to describe the server based distributed optimization like federated learning.
 Secondly, in the shared
message $v_{ij}^k$, the optimization variable $x_j^k$ and the
gradient $g_j^k$ are blended  together,  which makes it impossible for
a receiving agent to uniquely determine $x_j^k$ or $g_j^k$ based on
received information. In fact, the transmission of $x_i^k-\lambda^kg_i^k$ in our scheme amounts to transforming $x_i^k$ to a different reference frame  (by adding an unknown displacement $\lambda^kg_i^k$), which avoids the receiver from inferring the value of  $x_i^k$ or gradient $g_i^k$. In contrast,  the conventional scheme  directly discloses the optimization variable $x_i^k$, which also makes $g_i^k$ inferrable by an adversary that has  access to all messages shared in the network.} Note that, because $v_{ij}^k$ has the same
dimension as $x_j^k$, the new algorithm does not increase the
communication overhead compared with  conventional decentralized
algorithms. This one-variable only information-sharing scheme is
important in many applications such as machine learning because in
these applications the dimension of the optimization variables  can
scale to hundreds of millions, which causes significant
communication overhead and even communication bottlenecks
\cite{tang2020communication}.

 \noindent\rule{0.49\textwidth}{0.5pt}
\noindent\textbf{Algorithm 1: Decentralized nonconvex optimization
algorithm}

\vspace{-0.2cm}\noindent\rule{0.49\textwidth}{0.5pt}
\begin{enumerate}
    \item[] Parameters: $W$, $\lambda^k$
    \item {\bf for  $k=1,2,\cdots$ do}
    \begin{enumerate}
        \item Every agent $j$ computes and sends to  agent
        $i\in\mathbb{N}_j$
                \begin{equation}\label{eq:v_ij_DS}
             v^k_{ij}\triangleq w_{ij}(x_j^k-\lambda^kg_j^k)
        \end{equation}
        \item After receiving $v_{ij}^k$ from all $j\in\mathbb{N}_i$, agent $i$ updates its state as follows:
        \begin{equation}\label{update-rule-01}
            x_i^{k+1}=\sum_{j\in\mathbb{N}_i}v_{ij}^k=\sum\nolimits_{j\in\mathbb{N}_i}w^k_{ij}(x_j^k-\lambda^kg_j^k)
        \end{equation}
        \vspace{-0.4cm}
        \item {\bf end}
    \end{enumerate}
\end{enumerate}
\vspace{-0.2cm}\rule{0.49\textwidth}{0.5pt}

 Although the new algorithm  provides inherent privacy
 protection by sharing the mixture of the optimization variable and the
 gradient, the achieved privacy may not be strong enough. Therefore,
 we propose to inject additional additive noise to the gradient to
 ensure rigorous differential privacy. More specifically, instead of
 letting agent $j$
 send  $v_{ij}^k=w_{ij}(x_j^k-\lambda^kg_j^k)$ to agent $i$, we let agent $j$
 send $v_{ij}^k=w_{ij}(x_j^k-\lambda^k(g_j^k+n_j^k))$ to agent $i$, where $n_j^k\in\mathbb{R}^d$ is a $d$ dimensional Gaussian noise with mean zero and covariance matrix $\sigma I_d$. The
 detailed implementation procedure for individual agents is
provided in Algorithm 2.

 \noindent\rule{0.49\textwidth}{0.5pt}
\noindent\textbf{Algorithm 2: Decentralized nonconvex optimization
algorithm with differential privacy}

\vspace{-0.2cm}\noindent\rule{0.49\textwidth}{0.5pt}
\begin{enumerate}
    \item[] Parameters: $W$, $\lambda^k$
    \item {\bf for  $k=1,2,\cdots$ do}
    \begin{enumerate}
        \item Every agent $j$ computes and sends to  agent
        $i\in\mathbb{N}_j$
                \begin{equation}\label{eq:v_ij_DS}
             v^k_{ij}\triangleq w_{ij}(x_j^k-\lambda^k(g_j^k+n_j^k))
        \end{equation}
        \item After receiving $v_{ij}^k$ from all $j\in\mathbb{N}_i$, agent $i$ updates its state as follows:
        \begin{equation}\label{eq:update-rule-02}
            x_i^{k+1}=\sum_{j\in\mathbb{N}_i}v_{ij}^k=\sum\nolimits_{j\in\mathbb{N}_i}w^k_{ij}(x_j^k-\lambda^k(g_j^k+n_j^k))
        \end{equation}
        \vspace{-0.4cm}
        \item {\bf end}
    \end{enumerate}
\end{enumerate}
\vspace{-0.2cm}\rule{0.49\textwidth}{0.5pt}

  We will prove that not
 only can the noise $n_j^k$ ensure rigorous differential privacy for
 the data samples and
 gradient of agent $j$, but it will also bring  differential privacy protection to
 the optimization variable $x_j^k$ before the algorithm converges. (Note that
 we do not need to protect the privacy of the final optimization
 variable after convergence because the objective of decentralized
 optimization is to let individual agents learn the same optimal
 optimization
variable.) { In fact, as   illustrated later in Sec. V, compared with the conventional decentralized algorithm, our algorithm architecture greatly facilitates differential privacy design. For example,  to protect $g_i^k$ with a designated privacy strength,
since the transmitted message is $x_i^k-\lambda^kg_i^k$, the sender  can easily calculate the amount of noise that it should  add  to  shared messages. In contrast, in the conventional distributed optimization framework, as the shared information is $x_i^k$, it is not directly clear how much noise should be added to $x_i^k$ to achieve a certain privacy strength for $g_i^k$.}

 What is more interesting is that the  injected additive
noise does
 not compromise  the provable convergence of the algorithm, but instead, it
 ensures  avoidance of undesired stationary points like maxima and
 saddle points, and hence enhances the accuracy of decentralized
 nonconvex optimization. This is significant because
 differential-privacy noise is known to   sacrifice algorithmic accuracy
 for privacy. In the following two sections, we will rigorously analyze
 the convergence and maximum/saddle avoidance of Algorithm 2 and characterize its
 privacy-preservation performance in Sec. \ref{sec:convergence} and
 Sec. \ref{se:privacy_DS}, respectively.

%
%

For the convenience of the convergence analysis, we augment the
individual-agent dynamics in (\ref{eq:update-rule-02}) and obtain
the network-level dynamics:

\begin{equation}\label{eq:algorithm}
x^{k+1}=(W\otimes I_d)(x^k- \lambda^k (g^k+N^k))
\end{equation}
where  $\lambda^k\geq 0$ denotes the stepsize   at iteration $k$,
and $x^k$, $g^k$, and $N^k$ denote the stacked optimization
variables, gradients, and noise respectively, i.e.,
\[
\begin{aligned}
x^k&=[(x_1^k)^T,(x_2^k)^T,\cdots
(x_m^k)^T]^T , \\
g^k&=[(g_1^k)^T,(g_2^k)^T,\cdots (g_m^k)^T]^T ,\\
N^k&=[(n_1^k)^T,(n_2^k)^T,\cdots (n_m^k)^T]^T.
\end{aligned}
\]
The symbol $\otimes$ denotes the Kronecker product.

\section{Convergence Analysis}\label{sec:convergence}

In this section, we first prove that even in the presence of
differential-privacy noise, all $x_i^k$ in Algorithm 2 will reach
consensus almost surely (Sec. \ref{se:consensus}). Then, we will
  prove in Sec.
\ref{Sec:maxima} and Sec. \ref{sec:saddle} that the
differential-privacy noise in Algorithm 2 guarantees   avoidance of,
respectively, local maxima and saddle points.

For the convenience of  analysis, we define the average vector $\bar
x^k$ as
\[
\bar x^k=\frac{\sum_{i=1}^m x_i^k}{m}
\]
Since $W$ is column stochastic, from (\ref{eq:algorithm}), we have
\begin{equation}\label{eq:bar_x^k}
\bar x^{k+1}=\bar x^k - \frac{\lambda^k}{m} \sum_{i=1}^m
(g_i^k+n_i^k)
\end{equation}

\subsection{Consensus of all $x_i^k$}\label{se:consensus}
 We first prove that all $x_i^k$ $(1\leq i\leq m)$ converge to the average $\bar{x}^k$ almost
 surely. 

\begin{Theorem 1}\label{theo:consensus}
Under Assumption \ref{assumption:gradient}, Assumption \ref{As:bounded_gradient}, and Assumption
\ref{assumption:W}, when the stepsize $\lambda^k$ is square
summable, i.e., $\sum_{k=0}^{\infty}(\lambda^k)^2<\infty$, we have
$\lim_{k\rightarrow\infty}\|x_i^k-\bar{x}^k\|=0$ almost surely for
all $i$.
\end{Theorem 1}
Proof:
Using (\ref{eq:update-rule-02}) and (\ref{eq:bar_x^k}), one obtains
\begin{equation}\label{eq:x_error}
x^{k+1}-\bar{x}^{k+1}\otimes{\bf
1}=\bar{W}x^k-\lambda^k\bar{W}(g^k+N^k)
\end{equation}
where $\bar{W}=(W-\frac{{\bf 1}{\bf 1}^T}{m})\otimes{I_d}$.

Since  $(W-\frac{{\bf 1}{\bf 1}^T}{m}){\bf 1}=0$, we have for any
element $\bar{x}^k[\ell]$ of $\bar{x}^k$ (where $1\leq \ell\leq d$)
that, $(W-\frac{{\bf 1}{\bf 1}^T}{m}){\bf 1}\bar{x}^k[\ell]=0$
holds, which further leads to
\begin{equation}\label{eq:W}
(W-\frac{{\bf 1}{\bf 1}^T}{m})\otimes{I_d}\cdot\bar{x}^k \otimes{\bf
1}=0
\end{equation}

Combining (\ref{eq:x_error}) and (\ref{eq:W})
 yields
\begin{equation}
x^{k+1}-\bar{x}^{k+1}\otimes{\bf
1}=\bar{W}(x^k-\bar{x}^{k}\otimes{\bf 1})-\lambda^k\bar{W}(g^k+N^k)
\end{equation}
Using the definition of  $\eta=\|W-\frac{{\bf 1}{\bf 1}^T}{m}\|$ in
Assumption \ref{assumption:W}, we obtain
\begin{equation}\label{eq:norm_inequality}
\|x^{k+1}-\bar{x}^{k+1}\otimes{\bf 1}\|\leq
\eta\|x^k-\bar{x}^{k}\otimes{\bf 1}\|+\eta\lambda^k \|g^k+N^k\|
\end{equation}
By taking squares on both sides and using the inequality
\[
2ab\leq \epsilon a^2+\epsilon^{-1}b^2
\]
which holds for any $a,\,b\in\mathbb{R}$ and $\epsilon>0$, we obtain
\begin{equation}
\begin{aligned}
\|x^{k+1}-\bar{x}^{k+1}\otimes{\bf 1}\|^2\leq&
\eta^2(1+\epsilon)\|x^k-\bar{x}^{k}\otimes{\bf
1}\|^2\\
&+\eta^2(\lambda^k)^2(1+\epsilon^{-1}) \|g^k+N^k\|^2
\end{aligned}
\end{equation}
i.e.,
\begin{equation}
\begin{aligned}
\|x^{k+1}-\bar{x}^{k+1}\otimes{\bf 1}\|^2\leq&
\|x^{k}-\bar{x}^{k}\otimes{\bf 1}\|^2\\&-
\left(1-\eta^2(1+\epsilon)\right)\| x^k-\bar{x}^{k}\otimes{\bf
1} \|^2\\
&+\eta^2(\lambda^k)^2(1+\epsilon^{-1}) \|g^k+N^k\|^2
\end{aligned}
\end{equation}

By setting $1+\epsilon=\frac{1}{\eta}$ which further leads to
$\eta^2(1+\epsilon^{-1})=\frac{\eta^2}{1-\eta}$, one yields
\begin{equation}\label{eq:iteration}
\begin{aligned}
\|x^{k+1}-\bar{x}^{k+1}\otimes{\bf 1}\|^2\leq&
\|x^{k}-\bar{x}^{k}\otimes{\bf 1}\|^2\\&- \left(1-\eta\right)\|
x^k-\bar{x}^{k}\otimes{\bf
1} \|^2\\
&+ (\lambda^k)^2\frac{\eta^2}{1-\eta} \|g^k+N^k\|^2
\end{aligned}
\end{equation}

 Summing the preceding inequality over $k=0,1,\cdots$ yields
\begin{equation}\label{eq:consensus}
\begin{aligned}
&\left(1-\eta\right)\sum_{k=0}^{\infty}\|x^{k}-\bar{x}^{k}\otimes{\bf
1}\|^2+{\|x^{\infty}-\bar{x}^{\infty}\otimes{\bf 1}\|^2}\\
&-\|x^{0}-\bar{x}^{0}\otimes{\bf 1}\|^2\\
&\qquad\leq
\sum_{k=0}^{\infty}(\lambda^k)^2\frac{\eta^2}{1-\eta} \|g^k+N^k\|^2\\
&\qquad\leq \sum_{k=0}^{\infty}(\lambda^k)^2\frac{\eta^2}{1-\eta}
\|g^k\|^2+\sum_{k=0}^{\infty}(\lambda^k)^2\frac{\eta^2}{1-\eta}
\|N^k\|^2
\end{aligned}
\end{equation}

According to Assumption  \ref{As:bounded_gradient}, $g^k$ is bounded.
Therefore, we have
$\sum_{k=0}^{\infty}(\lambda^k)^2\frac{\eta^2}{1-\eta}
\|g^k\|^2<\infty$ when $\lambda^k$ is square summable. For the
second term on the right hand side of (\ref{eq:consensus}),
according to \cite{hall2014martingale}, we have $\|N^k\|^2$ being
finite almost surely (note that almost surely finite is different
from almost surely bounded). Therefore, under square summable
$\lambda^k$, the second term on the right hand side of
(\ref{eq:consensus}) is finite almost surely. In summary, the right
hand side of (\ref{eq:consensus}) is finite almost surely, meaning
that $(1-\eta)\sum_{k=0}^{\infty}\|x^{k}-\bar{x}^{k}\otimes{\bf
1}\|^2$ is finite almost surely, which further implies that
$\lim_{k\rightarrow\infty}\|x_i^k-\bar{x}^k\|=0$ holds almost surely
for all $i$.
\qed
\begin{Remark 1}
 The theorem can also be proven using our proof technique in \cite{wang2022tailoring}.
\end{Remark 1}

{
\begin{Remark 1}Besides almost sure convergence, we can also prove that   $\lim_{k\rightarrow \infty}\mathbb{E}\left[\|x_i^k-\bar{x}^k\|^2\right]\rightarrow 0$, where $\mathbb{E}[\cdot]$ is taken with respect to the $\sigma$-field generated by the Gaussian noise sequence $\{N^k\}$. More specifically, in the derivation of Theorem \ref{theo:consensus}, (\ref{eq:iteration}) also implies
\begin{equation}\label{eq:mean_squre}
\begin{aligned}
&\mathbb{E}\left[\|x^{k+1}-\bar{x}^{k+1}\otimes{\bf 1}\|^2\right]\\
&\leq
\mathbb{E}\left[\|x^{k}-\bar{x}^{k}\otimes{\bf 1}\|^2\right] - \left(1-\eta\right)\mathbb{E}\left[\|
x^k-\bar{x}^{k}\otimes{\bf
1} \|^2\right]\\
&\quad + (\lambda^k)^2\frac{\eta^2}{1-\eta} \mathbb{E}\left[\|g^k+N^k\|^2\right]\\
&\leq
\mathbb{E}\left[\|x^{k}-\bar{x}^{k}\otimes{\bf 1}\|^2\right] - \left(1-\eta\right)\mathbb{E}\left[\|
x^k-\bar{x}^{k}\otimes{\bf
1} \|^2\right]\\
&\quad + (\lambda^k)^2\frac{\eta^2}{1-\eta} \mathbb{E}\left[2\|g^k\|^2+2\|N^k\|^2\right]\\
&\leq
\mathbb{E}\left[\|x^{k}-\bar{x}^{k}\otimes{\bf 1}\|^2\right] - \left(1-\eta\right)\mathbb{E}\left[\|
x^k-\bar{x}^{k}\otimes{\bf
1} \|^2\right]\\
&\quad + (\lambda^k)^2\frac{\eta^2}{1-\eta}2mG^2 +(\lambda^k)^2\frac{\eta^2}{1-\eta}  2m\sigma
\end{aligned}
\end{equation}
where we have made use of  Assumption \ref{As:bounded_gradient} and the fact that $n_i^k$ has covariance matrix $\sigma I_d$.

 Summing the preceding inequality over $k=0,1,\cdots$ yields
 \begin{equation}\label{eq:consensus_2}
\begin{aligned}
&\left(1-\eta\right)\sum_{k=0}^{\infty}\mathbb{E}\left[\|x^{k}-\bar{x}^{k}\otimes{\bf
1}\|^2\right]+\mathbb{E}\left[\|x^{\infty}-\bar{x}^{\infty}\otimes{\bf 1}\|^2\right]\\
&-\mathbb{E}\left[\|x^{0}-\bar{x}^{0}\otimes{\bf 1}\|^2\right]\\
&\qquad\leq \sum_{k=0}^{\infty}(\lambda^k)^2\frac{\eta^2}{1-\eta}
2mG^2+\sum_{k=0}^{\infty}(\lambda^k)^2\frac{\eta^2}{1-\eta}
2m\sigma
\end{aligned}
\end{equation}
\end{Remark 1}

Therefore, the right
hand side of (\ref{eq:consensus_2}) is finite, meaning
that $(1-\eta)\sum_{k=0}^{\infty}\mathbb{E}\left[\|x^{k}-\bar{x}^{k}\otimes{\bf
1}\|^2\right]$ is finite, which further implies that
$\lim_{k\rightarrow\infty}\mathbb{E}\left[\|x_i^k-\bar{x}^k\|^2\right]=0$ holds
for all $i$.
}

\subsection{Avoidance of local maxima}\label{Sec:maxima}

Theorem \ref{theo:consensus} states that all $x_i^k$ will converge
to each other almost surely. However, it is still unclear what the
convergence point is. In the centralized case, it has been shown
that additive noise can enable a stochastic-approximation based
optimization process to converge to a local minimum when there are
no saddle points:

\begin{Lemma 1}\label{Lemma_avoiding_maximum}\cite{nevelson1976stochastic}
For a stochastic approximation process
\[
x^{k+1}=x^k-a^{k}(\nabla f(x^k)+q(k,x^k)+w^{k}),
\]
if the following conditions are satisfied:
\begin{enumerate}
\item $f(x)$ satisfies
$\lim_{\|x\|\rightarrow\infty}f(x)\rightarrow\infty$ and $\|\nabla
f(x)\|<C$ for some $C$;
\item $\sum_{k=0}^{\infty}a^{k}q(k,x^k)<\infty$ holds
almost surely;
\item $w^k$ are independent random variables
satisfying
 $\mathbb{E} \{w_i^k\}=0$ and $\mathbb{E} \{(w_i^k)^2\}=\sigma I$ with $\sigma<\infty$ for each element $w_i^k$ of $w^k$;
 \item  $a^k$ is not summable but square summable,
\end{enumerate}
  then state $x^k$ converges almost surely to a
point of the union of saddle and minima or the boundary of the
union. (It will avoid the local maxima almost surely).
\end{Lemma 1}

Recently, the above result has been extended in
\cite{tatarenko2017non} to
 push-sum based distributed optimization under the assumption of no saddle
points. However, the push-sum based distributed optimization
approach in \cite{tatarenko2017non}  has to  share two variables in
every iteration (one   optimization variable and  an additional
gradient-tracking variable), which is undesirable when the dimension
of the optimization variable is high. In fact, in modern deep
learning applications, the dimensions of optimization variables can
scale to hundreds of millions, and hence information sharing could
create significant communication overhead and even communication
bottlenecks \cite{tang2020communication}. In this paper, we show
that the result in Lemma \ref{Lemma_avoiding_maximum} can be
extended to our decentralized optimization framework which only
shares one variable in every iteration:

\begin{Theorem 1}\label{Theorem_x^i_to_bar{x}}
Under Assumption \ref{assumption:gradient},  Assumption \ref{As:bounded_gradient},  and Assumption
\ref{assumption:W}, when $\lambda^k$ is non-increasing, is not
summable but square summable, i.e.,
$\sum_{k=0}^{\infty}\lambda^k=\infty$ and
$\sum_{k=0}^{\infty}(\lambda^k)^2<\infty$, then  all  states $x^k_i$
converge almost surely to the same point in the union of saddles and
minima or the boundary of the union. (It will avoid  local maxima
almost surely).
\end{Theorem 1}

Proof:
It can be seen that under the conditions of the theorem, the
conditions in Theorem \ref{theo:consensus} are satisfied and hence
all $x_i^k$ will converge to the average state $\bar{x}^k$ almost
surely. So we only need to prove that $\bar{x}^k$ will converge to a
point in the union of saddles and minima or the boundary of the
union.

From  (\ref{eq:bar_x^k}), we can rewrite the dynamics of $\bar{x}^k$
as follows

\begin{equation}\label{eq:bar_x^k_with_q}
\begin{aligned}
\bar x^{k+1}&=\bar x^k -   \frac{\lambda^k}{m} \sum_{i=1}^m
(g_i^k+n_i^k)\\
 &=\bar x^k - \lambda^k\nabla
f(\bar{x}^k) -\lambda^k\frac{\sum_{i=1}^mn_i^k
}{m}\\
&\qquad +\lambda^k\left(\nabla f(\bar{x}^k)- \frac{\sum_{i=1}^mg_i^k
}{m} \right)
\end{aligned}
\end{equation}

Since when $n_i^k$ follows Gaussian distribution,
$\frac{\sum_{i=1}^mn_i^k }{m}$ also follows Gaussian distribution,
it can be seen that (\ref{eq:bar_x^k_with_q}) resembles the dynamics
in Lemma \ref{Lemma_avoiding_maximum} with
$q(k,x^k)=\frac{\sum_{i=1}^mg_i^k}{m}-\nabla f(\bar x^k)$ and
$a^k=\lambda^k$. Therefore, according to Lemma
\ref{Lemma_avoiding_maximum}, if we can prove that
$\sum_{k=0}^{\infty}\lambda^{k}q(k,x^k)$ is finite almost surely,
then it will follow that $\bar{x}^k$ will converge to a point in the
union of saddles and minima or the boundary of the union almost
surely, and hence   that all $x_i^k$ will converge to the same point
in the union of saddles and minima or the boundary of the union,
almost surely.

One can verify that $q(k,x^k)$ satisfies the following relationship:
\begin{equation}\label{eq:q_norm}
\begin{aligned}
\|q(k,x^k)\|&=\left\|\frac{\sum_{i=1}^mg_i^k}{m}-\nabla
f(\bar x^k)\right\|\\
&=\left\|\frac{\sum_{i=1}^m(g_i^k-\nabla f_i(\bar
x^k))}{m}\right\|\\
&\leq \frac{L}{m} \sum_{i=1}^m\|x_i^k-\bar{x}^k\| \\
&\leq \frac{L}{\sqrt{m}}\|x^k-\bar{x}^k\otimes {\bf 1}\|
\end{aligned}
\end{equation}
where the second to last inequality used the Lipschitz continuous
assumption of the gradients in Assumption \ref{assumption:gradient}
and the last inequality used the inequality $\sum_{i=1}^m a_i\leq
 \sqrt{m\sum_i^m (a_i)^2}$.

From (\ref{eq:norm_inequality}), we have
 \begin{equation}\label{eq:norm_inequality2}
 \begin{aligned}
\|&x^{k+1}-\bar{x}^{k+1}\otimes{\bf 1}\|\\
&\leq
\eta\|x^k-\bar{x}^{k}\otimes{\bf 1}\|+\eta\lambda^k \|g^k+N^k\|\\
&\leq \eta^2\|x^{k-1}-\bar{x}^{k-1}\otimes{\bf 1}\|+\eta^2\lambda^{k-1} \|g^{k-1}+N^{k-1}\|\\
&\qquad+\eta\lambda^k \|g^k+N^k\|\\
& \quad \vdots\\
&\leq \eta^{k+1}\|x^{0}-\bar{x}^{0}\otimes{\bf
1}\|+\sum_{l=0}^{k}\eta^{k+1-l} \lambda^{l} \|g^{l}+N^{l}\|
\end{aligned}
\end{equation}
which leads to
\begin{equation}\label{eq:sum_q}
\begin{aligned}
\sum_{k=0}^{\infty}\lambda^{k}q(k,x^k)\leq&
\frac{L}{\sqrt{m}}\sum_{k=0}^{\infty}\eta ^{k}\lambda^{k}
\|x^{0}-\bar{x}^{0}\otimes{\bf
1}\|\\
&+\frac{L}{\sqrt{m}}\sum_{k=0}^{\infty}\lambda^{k}\sum_{l=0}^{k-1}
\eta^{k-l} \lambda^{l} \|g^{l}+N^{l}\|
\end{aligned}
\end{equation}
Since $\lambda^k$ is non-increasing, we have
\[
\sum_{k=0}^{\infty}\lambda^{k}\eta^{k}\leq
\lambda^{0}\sum_{k=0}^{\infty}\eta^{k}=\lambda^{0}\frac{1}{1-\eta}<
\infty
\]
which further means that the first item on the right hand side of
(\ref{eq:sum_q}) is finite.

For the second term on the right hand side of (\ref{eq:sum_q}),
since $\{\lambda^k\}$ is a non-increasing sequence,  we always have
$\lambda^{k}\leq \lambda^{\ell}$ for $\ell\leq k$ and hence
\[
\sum_{k=0}^{\infty}\lambda^{k}\sum_{l=0}^{k-1} \eta^{k-l}
\lambda^{l}\leq \sum_{k=0}^{\infty} \sum_{l=0}^{k-1}\eta^{k-l}
(\lambda^{l})^2
\]
Noticing that $(\lambda^{l})^2$ is summable and $\eta$ resides in
the interval $(0,\,1)$, we have that
$\sum_{k=0}^{\infty}\lambda^{k}\sum_{l=0}^{k-1}\eta^{k-l}
\lambda^{l}$ is finite according to Lemma \ref{Lemma_Sundhar} in the
Appendix. Further using Assumption
\ref{As:bounded_gradient} that $g^k$ is always bounded and the
observation that $N^k$ is finite almost surely
\cite{hall2014martingale}, we have that the right hand side of
(\ref{eq:sum_q}) is finite almost surely. Therefore, we can conclude
that $\sum_{k=0}^{\infty}\lambda^{k}q(x^k,\bar{x}^k)$ is finite
almost surely.

In summary, under the conditions of the theorem, all conditions in
Lemma \ref{Lemma_avoiding_maximum} are satisfied. Therefore,  we can
conclude that $\bar{x}^k$ will converge to   a point in the union of
saddles and minima or the boundary of the union, and hence   all
$x_i^k$ will converge to the same point in the union of saddles and
minima or the boundary of the union.
\qed

\subsection{Avoidance of saddle points}\label{sec:saddle}

As we discussed in Sec. I, avoiding saddle points is a central
challenge for first-order based optimization methods. Recently some
advances have been reported on avoiding saddles in nonconvex
optimization (see e.g., \cite{ge2015escaping,jin2021nonconvex}).
However, these results are all for centralized optimization. Given
that in decentralized optimization   generally the local objective
functions of individual agents may have  saddle points different
from those of the aggregated objective function, and inter-agent
interactions also complicate the evolution of local optimization
variables, it is unclear if the results for the centralized case can
immediately be generalized to the
  decentralized optimization problem. Therefore, in this
section, we systematically address the saddle avoidance problem by
leveraging reported results on centralized optimization.

Inspired by the results in \cite{jin2021nonconvex}, we also use
coupling sequence to address the problem of saddle escaping:
\begin{Definition 1}
Given an optimization algorithm
\begin{equation}\label{eq:centralized_algorithm}
 x^{k+1}=x^k-\lambda^{k}(\nabla
f(x^k)+\xi^k(x^k)+w^{k})
\end{equation}
where $w^{k}$ is Gaussian noise with an identity covariance matrix
 and $\xi^k(x^k)$ is some function of the state $x^k$, we call $\{x'^k\}$
and $\{{x''^k}\}$
  coupling sequences starting from a strict saddle point $x_0$ if the following three conditions are
  satisfied:
\begin{enumerate}
\item   both sequences start from $x_0$;
\item both are obtained as separate runs
  of the optimization algorithm under the same $\xi^k(x^k)$;
\item  both are obtained under $w'^{k}$ and $w''^k$ that are only different in the $e_1^T$ direction,
  i.e.,
  $e_1^Tw'^k=-e_1^Tw''^k$, where $e_1$ denotes the eigenvector associated with the minimum eigenvalue of the Hessian  matrix at the saddle point $x_0$.
\end{enumerate}
\end{Definition 1}

The results for the centralized optimization in
\cite{jin2021nonconvex} show that for a strict saddle $x_0$, if with
a positive probability, the magnitude of the projected $\xi^k(x^k)$
on $x'^k-x''^k$ is less than half of the magnitude of the projected
  $w^k$ on $x'^k-x''^k$, then  the algorithm in (\ref{eq:centralized_algorithm}) can
effectively avoid the saddle point:

\begin{Lemma 1}\label{lemma_saddle_avoiding}\cite{jin2021nonconvex}
For the stochastic approximation process $x^{k+1}=x^k-a (\nabla
f(x^k)+\xi^k(x^k)+w^{k})$ where $w^{k}$ is  Gaussian noise with an
identity covariance matrix, suppose that  $\{x'^k\}$ and
$\{{x''^k}\}$ are
  coupling sequences starting from a strict saddle point $x_0$. If
  with a positive probability,
  the magnitude of the projected   $\xi^k(x^k)$ on   $x'^k-x''^k$ is
less than half of the magnitude of projected   $w^k$ on $x'^k-x''^k$
(the dynamics of the difference $x'^k-x''^k$ is dominated by $w^k$),
then for any given probability $0<\mu<1$, $x^k$ will escape the
saddle $x_0$ with probability at least $1-\mu$ after at most
$\mathcal{O}(\frac{\log(\frac{1}{\mu})}{a})$ iterates with a
sufficiently small constant stepsize $a \leq \frac{1}{\ell}$.
\end{Lemma 1}

According to Theorem \ref{theo:consensus}, under Assumption
\ref{assumption:W}, all $x_i^k$ will converge to the mean state
$\bar{x}^k$ almost surely, so we only need to consider if
$\bar{x}^k$ can avoid saddle points. Given that the dynamics of
$\bar{x}^k$ is governed by (\ref{eq:bar_x^k_with_q}), from Lemma
\ref{lemma_saddle_avoiding}, if we can prove that   the difference
between two coupling sequences initiating from a saddle point is
governed by $\frac{\sum_{i=1}^mn_i^k }{m}$ rather than
$q(k,x^k)=\nabla f(\bar{x}^k)- \frac{\sum_{i=1}^mg_i^k }{m} $, then
$\bar{x}^k$ will escape the saddle point, meaning that all $x_i^k$
will escape the saddle point. More specifically, we can prove the
following result:
\begin{Theorem 2}\label{theo:saddle}
Under Assumption \ref{assumption:gradient}, Assumption \ref{As:bounded_gradient}, and Assumption
\ref{assumption:W}, the differentially-private decentralized nonconvex optimization
algorithm in Algorithm 2 avoids saddle points with probability at
least $1-\mu$ for any $0<\mu<1$ under the  following stepsize
strategy:
\begin{enumerate}
\item constant and small
enough $\lambda^k \leq \frac{1}{\ell}$  for the first
$\mathcal{O}(\frac{\log(\frac{1}{\mu})}{\lambda})$ iterates; and
then
\item diminishing $\lambda^k$  satisfying the non-summable but square summable
condition.
\end{enumerate}
\end{Theorem 2}

Proof:
We first study the dynamics of the difference between two coupling
sequences $\bar{x}'^k$ and $\bar{x}''^k$. Since the evolution of
$\bar{x}^k$ is governed by
\[
\bar x^{k+1}=\bar x^k  - \lambda^k\nabla f(\bar{x}^k) +\lambda^k
q(k,x^k)-\lambda^k\frac{\sum_{i=1}^mn_i^k }{m}
\]
we can represent the dynamics of the two coupling sequences ${x}'^k$
and  ${x}''^k$ as
\[
\begin{aligned}
 x'^{k+1}&=x'^k  - \lambda^k\nabla f({x}'^k) +\lambda^k
q(k,x'^k)-\lambda^k\xi'^k\\
 x''^{k+1}&=x''^k  - \lambda^k\nabla f({x}''^k) +\lambda^k
q(k,x''^k)-\lambda^k\xi''^k
\end{aligned}
\]
where $\xi'^k$ and $\xi''^k$ represent the corresponding aggregated
Gaussian noise  $\frac{\sum_{i=1}^mn_i^k }{m}$ differing in only the
$e_1$ direction.

Therefore, the difference $\tilde{x}^k\triangleq\xi'^k-\xi''^k$ has
the following dynamics
\begin{equation}\label{eq:difference_dynamics}
\begin{aligned}
\tilde{x}^{k+1}=&\tilde{x}^{k}-\lambda^k(\nabla f({x}'^k)-\nabla
f({x}''^k))\\
&+\lambda^k(q(k,x'^k)-q(k,x''^k))-\lambda^k(\xi'^k-\xi''^k)
\end{aligned}
\end{equation}
Let  $\tilde{\xi}^k\triangleq\xi'^k-\xi''^k$ and
$\Delta^k\triangleq\int_{0}^1\nabla^2f(\phi
x'^k-(1-\phi)x''^k)d\phi$. Noticing that $\nabla f({x}'^k)-\nabla
f({x}''^k)=\int_{0}^1\nabla^2f(\phi
x'^k-(1-\phi)x''^k)d\phi(x'^k-x''^k)$, we can rewrite
(\ref{eq:difference_dynamics}) as
\[
\tilde{x}^{k+1}=(I-\lambda^k\Delta^k)
\tilde{x}^{k}+\lambda^k(q(k,x'^k)-q(k,x''^k))-\lambda^k\tilde{\xi}^k
\]

From (\ref{eq:q_norm}) and (\ref{eq:norm_inequality2}), we have
\[
\begin{aligned}
\|q(k,x'^k)\| &\leq \frac{L}{\sqrt{m}}\|x^k-\bar{x}^k\otimes {\bf
1}\|\\
&\leq \frac{L\eta^{k}}{\sqrt{m}}\|x^{0}-\bar{x}^{0}\otimes{\bf 1}\|\\
&+
\frac{L }{\sqrt{m}}\sum_{l=0}^{k} \eta^{k-l}\lambda^{l}
\|g^{l}+N^{l}\|
\end{aligned}
\]
Since $g^l$ is   bounded according to Assumption
\ref{As:bounded_gradient} (denote the upper bound as $G$), and
under a positive probability, $N^l$ is less than some positive $T$,
with a positive probability, $\|g^{l}+N^{l}\|$ is less than $G+T$,
which means that we have the following relationship with a positive
probability under a constant stepsize $\lambda$:
\[
\|q(k,x'^k)\| \leq
\frac{L\eta^{k}}{\sqrt{m}}\|x^{0}-\bar{x}^{0}\otimes{\bf 1}\|+
\frac{L\lambda  }{\sqrt{m}}\frac{1-\eta^k}{1-\eta} (G+T)
\]

Similarly, we   have
\[
\|q(k,x''^k)\| \leq
\frac{L\eta^{k}}{\sqrt{m}}\|x^{0}-\bar{x}^{0}\otimes{\bf 1}\|+
\frac{L\lambda  }{\sqrt{m}}\frac{1-\eta^k}{1-\eta} (G+T)
\]
with a positive probability. Therefore, with a positive probability,
we have
\[
\begin{aligned}
\|q(k,x'^k)-q(k,x''^k)\| \leq&
\frac{2L\eta^{k}}{\sqrt{m}}\|x^{0}-\bar{x}^{0}\otimes{\bf 1}\|\\
&+ \frac{2L\lambda  }{\sqrt{m}}\frac{1-\eta^k}{1-\eta} (G+T)
\end{aligned}
\]
In the mean time, given that $\tilde{\xi}^k$ is also Gaussian, we
have that with a positive probability,
\[
\|\tilde{\xi}^k\|>T
\]
holds.

Therefore, for $\lambda$ sufficiently small, we can always have
$\|\tilde{\xi}^k\|>2\|q(k,x'^k)-q(k,x''^k)\|$ with a positive
probability, implying that with a positive probability,
$\tilde{\xi}^k$ dominates the dynamics of $x'^k-x''^k$. Then
according to Lemma \ref{lemma_saddle_avoiding}, we have that
$\bar{x}^k$  can avoid the saddle with at least probability $1-\mu$
for any $0<\mu<1$. Further using the result from Theorem
\ref{theo:consensus}  that all $x_i^k$ converge to $\bar{x}^k$
almost surely yields that all $x_i^k$ can avoid the saddle with at
least probability $1-\mu$ for any $0<\mu<1$.
\qed


\section{Privacy Analysis}\label{se:privacy_DS}


In this section, we prove that Algorithm 2 can provide rigorous
differential privacy for data samples
 $s_{ij}$, individual agents' gradients $g_i^k$,  and intermediate optimization
variables $x_i^k$. By intermediate optimization variables, we mean
the evolution of optimization variables  $x_i^k$ before convergence
is achieved. Note that after convergence, it is the exact objective
of decentralized optimization  to have all agents arrive at the same
optimization variable (and hence know each other's value).
 We first analyze the achieved privacy strength for data samples.

According to differential privacy, the minimum of noise variance
required to achieve $(\epsilon,\delta)$-differential privacy for
data samples is determined by the  sensitivity function
\begin{equation}\label{eq:sentivity}
S_{s,i}=\sup\limits_{\|s_{i,p}-s_{i,q}\|_1\leq 1}\|
M_{i}^k(s_{i,p})-M_i^k(s_{i,q})\|_1
\end{equation}
where $M_i^k\triangleq x_i^k-\lambda^kg_i^k $. Note that we replaced
$v_{ji}^k=w_{ji}(x_i^k-\lambda^kg_i^k)$ with $M_i^k$ to calculate
the sensitivity function because the coefficient $w_{ji}$ is
publicly known. Further notice that $g_i^k=\nabla f_i(x_i^k)$ holds
and in the above sensitivity function, changing one data sample from
$s_{i,p}$ to $s_{i,q}$ only affects one cost function
$\ell_i(\cdot,\cdot)$, and thus according to Assumption
\ref{assumption:gradient} and the relationship between
$\ell_i(\cdot,\cdot)$ and $f_i(\cdot)$ in (\ref{eq:f_i}), one can
obtain that $S_i=\frac{\nu\lambda^k}{n_i}$. Then making use of the
standard result  in differential privacy, we have the following
theorem for privacy protection on data samples:

\begin{Theorem 1}\label{theo:differential_privacy}
For any $\epsilon,\delta\in(0,\,1)$, at iteration $k$, the noise
$n_i^k$ can ensure $(\epsilon,\delta)$-differential privacy for
agent $i$'s every data sample $s_{i,q}$  when the variance
$\sigma^2$ satisfies
\begin{equation}\label{eq:differtial_privacy}
\sigma^2\geq
2\nu^2(\lambda^k)^2\frac{\ln(1.25/\delta)}{n_i^2\epsilon^2}
\end{equation}
\end{Theorem 1}\label{th:differential_privacy}
Proof:
According to \cite{dwork2014algorithmic}, a Gaussian noise of
variance  $\sigma^2\geq 2\frac{\ln(1.25/\delta)(S_f)^2}{\epsilon^2}$
can achieve $(\epsilon,\delta)$-differential privacy  for any
$\epsilon,\delta\in(0,\,1)$ where $S_f$ denotes the sensitivity
function. Thus the proof will be completed by  incorporating the
sensitivity function value of $S_i=\frac{\nu\lambda_i^k}{n_i}$
obtained just above the statement of the theorem.
\qed
\begin{Remark 1}
Note that there are infinitely many $(\epsilon,\,\delta)$ pairs that
satisfy (\ref{eq:differtial_privacy}) in Theorem
\ref{theo:differential_privacy}.
\end{Remark 1}

Note that not only does the   differential-privacy noise $n_i^k$
added to agent $i$  enable  privacy protection for agent $i$'s data
samples, but it also provides differential-privacy protection for
agent $i$'s gradient. The sensitivity function for agent $i$'s
gradient is given by
\begin{equation}\label{eq:sentivity}
S_{g,i}=\sup\limits_{\|g_{i}^k-{g'}_{i}^k\|_1\leq
1}\|M_{i}^k(g_{i}^k)-M_i^k({g'}_{i}^k)\|_1
\end{equation}
where $M_i^k\triangleq x_i^k-\lambda^kg_i^k$.

It can be shown that $S_{g,i}=\lambda^k$, and hence we  have the
following theorem:
\begin{Theorem 1}\label{theo:differential_privacy_gradient}
For any $\epsilon,\delta\in(0,\,1)$, at iteration $k$, the noise
$n_i^k$ can also ensure $(\epsilon,\delta)$-differential privacy for
agent $i$'s gradient $g_{i}^k$  when the variance $\sigma^2$
satisfies
\[
\sigma^2\geq 2(\lambda^k)^2\frac{\ln(1.25/\delta)}{\epsilon^2}
\]
\end{Theorem 1}\label{th:differential_privacy}
Proof:
The result follows from a similar line of argument as in the proof
of  Theorem \ref{theo:differential_privacy}.
\qed

 From Theorem \ref{theo:differential_privacy}, we can see that  to achieve a fixed level of
differential privacy for data samples, the required noise level
decreases with an increase in the number of data samples   $n_i$. To
the contrary, Theorem \ref{theo:differential_privacy_gradient} shows
that the required noise level for differential privacy of gradients
is not affected by the number of samples, which is understandable
since the gradient of an agent is always computed from all data
samples of the agent in gradient descent algorithms.

The same noise $n_i^k$ also provides privacy protection for
optimization variables. In fact, we can see that the amount of noise
applied on $x_i^k$ is $\lambda^kn_i^k$   in shared information. And
the variance of this noise is $(\lambda^k)^2\sigma^2$. Since it can
be verified that   the  sensitivity function for $x_i^k$ is
\begin{equation}\label{eq:sentivity}
S_{x,i}=\sup\limits_{\|x_{i}^k-{x'}_{i}^k\|_1\leq
1}\|M_{i}^k(x_{i}^k)-M_i^k({x'}_{i}^k)\|_1=1
\end{equation}
where $M_i^k\triangleq x_i^k-\lambda^kg_i^k$,  we can obtain that
the same noise $n_i^k$ also enables the following differential
privacy for individual agents' optimization variables:

\begin{Theorem 1}
At iteration $k$, the noise $n_i^k$ also ensures
$(\epsilon,\delta)$-differential privacy for agent $i$'s
optimization variable $x_i^k$ for any   $\epsilon,\delta\in(0,\,1)$
when the variance $\sigma^2$ satisfies $\sigma^2\geq
2\frac{\ln(1.25/\delta)}{(\lambda^k)^2\epsilon^2}$.
\end{Theorem 1}
Proof:
The result follows from a similar line of argument as in the proof
of Theorem \ref{theo:differential_privacy}.
\qed

\begin{Remark 1}
Note that different from the enabled privacy for data samples and
gradients, under a fixed noise variance $\sigma$, the strength of
enabled $(\epsilon,\delta)$-differential privacy for $x_i^k$
decreases with a decrease in $\lambda^k$. When $\lambda^k$ tends to
zero, the strength of enabled privacy for $x_i^k$ will decrease to
zero. However, note that this is acceptable since the purpose of
decentralized optimization is for all agents to learn the same
optimum value for the optimization variable   cooperatively.
\end{Remark 1}

{
\begin{Remark 1}
  Also note that in all the above results on sensitivity and differential privacy, we have assumed that the adversary knows the underlying algorithm and can observe every shared message in the network, namely, the adversary can launch both the honest-but-curious attack and the eavesdropping attack discussed in Sec. \ref{se:privacy_formulation}. If the adversary is weaker in the sense that it  can only launch the honest-cut-curious attack (can only observe messages shared on some but not all links), then the sensitivity of agents whose messages are unaccessible to the adversary will be smaller (and even zero), and hence these agents will have a stronger privacy protection against the adversary. The same conclusion can be drawn  for the case where the adversary  can only launch the eavesdropping attack  (does not know the underlying algorithm).
\end{Remark 1}
}

\section{Numerical Experiments}

In this section, we evaluate the performance of the proposed
decentralized  optimization algorithm using numerical experiments in
decentralized nonconvex optimization applications. More
specifically, we evaluate the performance of the proposed algorithm
using two application scenarios, one in decentralized estimation and
the other in Independent Component Analysis, a popular dimension
reduction tool in statistical machine learning and signal processing
\cite{hyvarinen2000independent}.

\subsection{Decentralized estimation based numerical experiments}

We consider a canonical decentralized estimation problem where a
 network of $m$ sensors collectively estimate an  unknown parameter
$\theta\in\mathbb{R}^d$. More specifically, we assume that each
sensor $i$ has a measurement  of the parameter,
$Y_{i}=M\theta+w_{i}$, where $M\in\mathbb{R}^{s\times d}$ is the
measurement matrix of agent $i$ and $w_{i}$ is measurement noise.
Then the estimation of parameter $\theta$ can be solved using the
optimization problem formulated as
(\ref{eq:optimization_formulation1}), with each $f_i(\theta)$ given
as
\[
f_i(\theta)=\|Y_i-M\theta\|^2+\kappa\|\theta\|^3
\]
Here $\kappa$ is a  regularization parameter, which will be chosen
to have some desired properties for $f_i(\cdot)$.

It can be verified that when $\kappa$ is a positive number,
$f_i(\cdot)$ will be a convex function. In order to have a nonconvex
objective function so as to test and evaluate the performance of our
algorithm in nonconvex optimization, we set $\kappa$ as
$\kappa=-0.1$. We took $M$ and $Y_i$ ($1\leq i\leq 5$) as follows
\[
M=\left[\begin{array}{cc}1&0\\0&2\\0&0\end{array}\right],\quad
Y_i=i\times \left[\begin{array}{c} 1/3\\2/3\\0\end{array}\right]
\]

Using the relationship
\[
\nabla f_i=-2M^TY_i+2M^TM\theta+3\kappa\|\theta\|\theta
\]
and
\[
\nabla^2
f_i=2M^TM-3\kappa\|\theta\|I_2+3\kappa\frac{1}{\|\theta\|}\theta\theta^T
\]
one can obtain that the aggregated function
$f(\theta)=\sum_{i=1}^5\frac{ f_i(\theta)}{5}$ has a local minimum
  at $\theta=\left[\begin{array}{c}1.3478,
\\1.0690 \end{array}\right]$ and a
saddle point at $\theta =\left[\begin{array}{c} -7.4336\\
1.3959\end{array}\right]$.

\begin{figure}
\hspace{-0.3cm}
\includegraphics[width=0.5\textwidth]{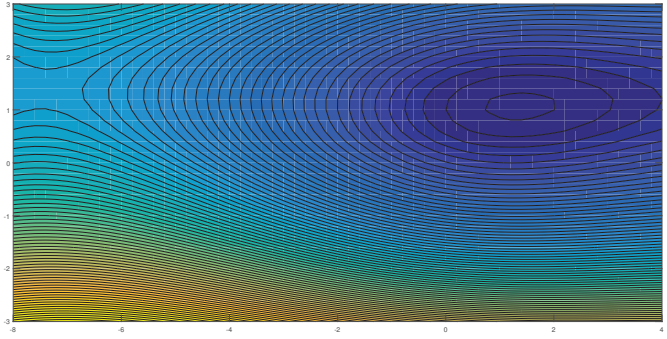}
    \caption{ A two-dimensional contour graph  of $f(\theta)$ in $[-8,\,4]\times [-3,\,3]$.}
    \label{fig:countour}
\end{figure}

To facilitate numerical experiments, we focus on the region
$\theta\in[-8,\,4]\times [-3,\,3]$. Outside this region
 we manipulate $f_i(\theta)$ to make it increase  linearly
with $\|\theta\|$ with continuous and smooth connection on the
boundary of $[-8,\,4]\times [-3,\,3]$. By doing so, our optimization
problem has one minimum at $\theta=\left[\begin{array}{c}1.3478,
\\1.0690 \end{array}\right]$ (and no other local minimum) and one saddle point at   $\theta=\left[\begin{array}{c} -7.4336\\
1.3959\end{array}\right]$, and it can be verified that the saddle
point is a strict saddle point. Please see Fig. \ref{fig:countour}
for a two-dimensional contour graph of $f(\theta)$ on
$[-8,\,4]\times [-3,\,3]$.

\begin{figure}
    \begin{center}
        \includegraphics[width=0.25\textwidth]{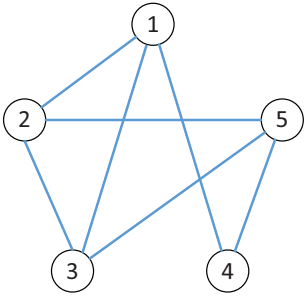}
    \end{center}
    \caption{The interaction topology of the network.}
    \label{fig:topology}
\end{figure}

\begin{figure}
    \begin{center}
        \includegraphics[width=0.5\textwidth]{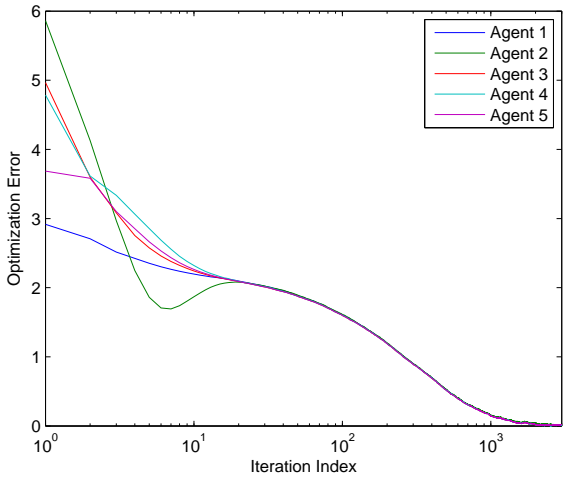}
    \end{center}
    \caption{Evolution of all agents' optimization errors  when initialized with random  values.}
    \label{fig:random_initialization}
\end{figure}

\begin{figure}
    \begin{center}
        \includegraphics[width=0.5\textwidth]{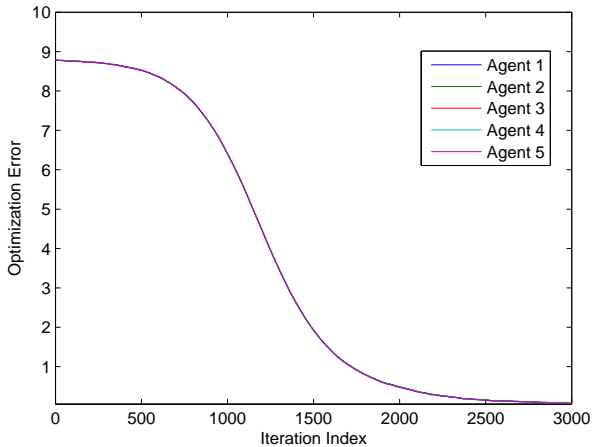}
    \end{center}
    \caption{Evolution of all agents' optimization errors  when  every agent was initialized from the saddle point.}
    \label{fig:initialization_saddle}
\end{figure}

We considered a network of five agents  interacting on a  graph
depicted in Fig. \ref{fig:topology}. In the numerical experiments,
we set the stepsize as $\lambda^k=0.02$ for $k\leq 500$ and switched
it to $\lambda_i=\frac{1}{k}$ for $k>500$. $w_{i}$ was  uniformly
distributed in $[0,\,1]$. We added Gaussian noise with zero mean and
variance $\sigma=0.5$ in the gradients  for the purpose of both
privacy protection and global convergence. We first initialized the
optimization variables  randomly to check if the algorithm can
guarantee consensus in decentralized optimization in the presence of
differential-privacy noise. The evolution of all agents'
optimization variables is illustrated in Fig.
\ref{fig:random_initialization}, which confirms that the algorithm
can indeed ensure all agents to converge to the same state under
differential-privacy noise. To evaluate the performance of our
algorithm with regard to  avoiding the saddle point,
we also initialized all agents on the saddle point, i.e.,  $x_i^0=\left[\begin{array}{c} -7.4336\\
1.3959\end{array}\right]$ for all $1\leq i\leq 5$. Clearly, without
the differential-privacy noise, all states would be trapped at the
saddle point. In contrast,  the differential-privacy noise
avoided the saddle point and ensured the convergence of all agents
to the optimal value, as illustrated in Fig.
\ref{fig:initialization_saddle}, which corroborates the theoretical
results in Theorem \ref{theo:saddle}. {We have also evaluated the influence of the magnitude of $\sigma$ on the optimization error. The results are summarized in Table I, where each data point was the average of 100 runs. It can be seen that the optimization error increases with an increase in the noise magnitude $\sigma$.
}

\begin{table*}
\centering
\caption{Final optimization error under different $\sigma$ at $k=3000$}
\begin{tabular}{  c|c|c|c|c|c|c  }
  &$\sigma=0.1$& $\sigma=0.2$ & $\sigma=0.3$&$\sigma=0.4$& $\sigma=0.5$&$\sigma=0.6$ \\
  \hline
  Average optimization error&0.048&0.058&0.064& 0.070 & 0.078 &0.091
\end{tabular}
\end{table*}

\subsection{Independent Component Analysis based numerical experiments}
Independent Component Analysis (ICA) is   widely used in signal
processing and statistical machine learning to reduce the dimension
of data \cite{hyvarinen2000independent}. Modeling  the data vector
as $Y=AZ$ with $A\in\mathbb{R}^{d\times d}$ an orthonormal matrix
and $Z\in\mathbb{R}^{d}$ a non-Gaussian random sample  of $d$
independent entries, the objective of ICA is to recover one of
multiple columns of $A$ from independent observations
$Y\in\mathbb{R}^{d}$. In standard practice, the  $Y$ vector should
be whitened to have zero mean and an identity covariance matrix.
Furthermore, in standard practice,  the elements of $Z$ are usually
assumed to have a fourth moment $\mu\neq 3$, and   the $d$ columns
of $A$ are usually denoted as ${\bf a}_1,\,{\bf a}_2\,\cdots,{\bf
a}_d$.
 Under these conditions, ICA is usually cast as the following
optimization problem:
\[
\min\limits_{\|u\|=1} -{\rm sign}(\mu-3)\cdot
\mathbb{E}\left[(u^TY)^4\right
]
\]

\begin{figure}
    \begin{center}
        \includegraphics[width=0.5\textwidth]{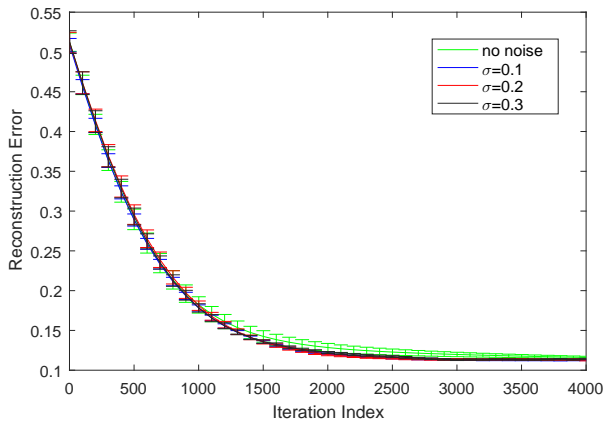}
    \end{center}
    \caption{Evolution of the maximum reconstruction error among all agents in the ICA application.}
    \label{fig:ICA}
\end{figure}

The saddles of the above optimization problem include (but not
limited to) all $u^{\ast}=d^{-1/2}(\pm 1,\cdots,\pm 1)$, all of
which satisfy the strict-saddle condition \cite{li2016online}. We
implemented our algorithm to solve the above optimization problem.
(It is worth noting that the equality constraint $\|u\|=1$ can be
handled by using the method of Lagrange multipliers
\cite{ge2015escaping}.) In the implementation, we set $d=10$. We
also generated 800 random samples $Y$ by randomly selecting each
entry of $Z$ from a uniform distribution in $\{-1,\,1\}$. The 800
samples were evenly distributed among the five agents (each agent
had 160 samples). We set the stepsize as $\lambda^k=0.003$ for
$k\leq 100$ and switched it to $\lambda_i=\frac{3}{10k}$ for
$k>100$. The network interaction topology is still the same as in
Fig. \ref{fig:topology}. We evaluated the performance of our
algorithm under different variances $\sigma$ of the
differential-privacy noise. In each case, we ran the algorithm for
100 times. In each of the 100 runs, we randomly set $A$ to a  random
orthonormal matrix and randomly selected $x_i^0$. The evolution of
the maximal reconstruction error of the first column of A,
$\bf{a}_1$, among all five agents is shown in Fig. \ref{fig:ICA}. It
can be seen that after adding differential-privacy noise, our
algorithm obtains comparable or even slightly better convergence
accuracy and speed compared with the noise-free case. {This is understandable since our algorithm can always avoid saddle points, whereas the noise-free case was trapped at saddle points in some of the 100 runs. Note that since the initial conditions for the 100 runs were randomly selected, the noise-free case is not always trapped at saddle points. }

\section{Conclusions}

This paper has proposed an algorithm for  decentralized nonconvex
optimization  that can achieve rigorous differential privacy with
guaranteed convergence   to a minimum point.   By leveraging
diminishing stepsizes, the algorithm avoids sacrificing provable convergence for differential privacy, which is a common problem with
existing differential-privacy based algorithms for decentralized
optimization. Besides enabling privacy protection for data samples
and gradients, the approach  also achieves privacy protection for
optimization variables until the algorithm converges. More
interestingly, we have  proved that our algorithm has guaranteed
saddle avoidance   in a polylogrithmic number of iterations. The
guarantee on  differential privacy, algorithmic convergence, and
saddle-avoidance  simultaneously has not been reported in
decentralized optimization literature. Note that since in
decentralized optimization
 individual agents may have   saddle points different from those of  the centralized counterpart,
 the
saddle-avoidance result obtained for decentralized optimization is
highly nontrivial compared with existing results for centralized
optimization. Numerical experiments for both a decentralized
estimation problem and an independent component analysis problem
confirm the effectiveness of the proposed algorithm.

\section*{Appendix A}

\begin{Lemma 1}\label{Lemma_Sundhar}\cite{ram2010distributed}
Let $\{\gamma^k\}$ be a scalar sequence. If $\gamma^k\geq 0$ for all
$k$, $\sum_{k=0}^{\infty}\gamma^k<\infty$, and $0<\beta<1$, then
$\sum_{k=0}^{\infty}(\sum_{\ell=0}^k(\beta)^{k-\ell}\gamma^{\ell})<\infty$.
\end{Lemma 1}


\bibliographystyle{plain}        
\bibliography{reference1}           

\end{document}